\magnification=1200
\overfullrule=0pt
\centerline {\bf Another four critical points theorem}\par
\bigskip
\bigskip
\centerline {BIAGIO RICCERI}\par
\bigskip
\bigskip
\centerline {\it Dedicated to the memory of Ky Fan,
with my immense esteem and admiration}\par
\bigskip
\bigskip
{\bf Abstract:} In this paper, making use of Theorem 2 of [5], we
establish a new four critical points theorem which can be regarded as a
companion to Theorem 1 of [4]. We also present an
 application to the Dirichlet problem for a class of quasilinear
elliptic  equations.\par
\par
\bigskip
\bigskip
\bigskip
\bigskip
The aim of this paper is to establish a new four critical points
theorem (Theorem 1 below) that can be regarded as a companion to
Theorem 1 of [4].\par
\smallskip
As in [4], our key tool is the multiplicity result on global minima established
in [3]. The use of such a result requires the validity of a strict minimax inequality
which is explicitly assumed in [4].\par
\smallskip
A novelty of Theorem 1 is that no minimax inequality appears among the hypotheses.
This is possible thanks to the use of Theorem 2 of the very recent
[5] which just highlights a rather
general situation where the strict minimax inequality occurs.\par
\smallskip
In other words, Theorem 1 should be regarded as the fruit that one obtains by combining
the underlying ideas of [4] with Theorem 2 of [5].\par
\smallskip
For the reader convenience, we start just recalling Theorem 2 of [5].\par
\smallskip
First, we introduce the following notations.\par
\smallskip
If $X$ is a non-empty set and $\Gamma, \Psi, \Phi:X\to {\bf R}$ are
three given functions,
for each $\mu>0$ and $r\in ]\inf_X\Phi,\sup_X\Phi[$, we put\par
$$\alpha(\mu \Gamma+\Psi,\Phi,r)=\inf_{x\in \Phi^{-1}(]-\infty,r[)}
{{\mu \Gamma(x)+\Psi(x)-\inf_{\Phi^{-1}(]-\infty,r])}(\mu \Gamma+\Psi)}
\over {r-\Phi(x)}}$$
and
$$\beta(\mu \Gamma+\Psi,\Phi,r)=\sup_{x\in \Phi^{-1}(]r,+\infty[)}
{{\mu \Gamma(x)+\Psi(x)-\inf_{\Phi^{-1}(]-\infty,r])}
(\mu \Gamma+\Psi)}\over {r-\Phi(x)}}\ .$$
When $\Psi+\Phi$ is bounded below, for each $r\in ]\inf_X\Phi,\sup_X\Phi[$ such that
$$\inf_{x\in \Phi^{-1}(]-\infty,r])}\Gamma(x)<
\inf_{x\in \Phi^{-1}(r)}\Gamma(x)
\ ,$$
we put
$$\mu^*(\Gamma,\Psi,\Phi,r)=\inf\left \{ {{\Psi(x)-\gamma+r}\over
{\eta_r-\Gamma(x)}} :
x\in X, \Phi(x)<r, \Gamma(x)<\eta_r\right \}\ ,$$
where
 $$\gamma=\inf_{x\in X}(\Psi(x)+\Phi(x))$$
and
$$\eta_r=\inf_{x\in \Phi^{-1}(r)}\Gamma(x)
\ .$$

THEOREM A ([5], Theorem 2). - {\it Let $X$ be a topological space and $\Gamma, \Psi, \Phi:X\to
{\bf R}$ three sequentially lower
semicontinuous functions, with $\Gamma$ also sequentially inf-compact,
 satisfying the following conditions:\par
\noindent
$(a)$\hskip 5pt $\inf_{x\in X}(\mu \Gamma(x)+\Psi(x))=-\infty$ for all $\mu>0$\ ;
\par
\noindent
$(b)$\hskip 5pt $\inf_{x\in X}(\Psi(x)+\Phi(x))>-\infty$\ ;\par
\noindent
$(c)$\hskip 5pt there exists $r\in ]\inf_X \Phi,\sup_X \Phi[$ such that
$$\inf_{x\in \Phi^{-1}(]-\infty,r])}\Gamma(x)<
\inf_{x\in \Phi^{-1}(r)}\Gamma(x)\ .$$
Under such hypotheses, for each
$\mu>\max\{0,\mu^*(\Gamma,\Psi,\Phi,r)\}$, one has
$$\alpha(\mu \Gamma+\Psi,\Phi,r)=0$$
and
$$\beta(\mu \Gamma+\Psi,\Phi,r)>0\ .$$}
\medskip
As we said, the key tool in the proof of Theorem 1 is provided by the following
\medskip
THEOREM B ([3], Theorem 1). - {\it Let $X$ be a topological space,
$A\subseteq {\bf R}$ an open interval and $P:X\times A\to {\bf R}$ a
function satisfying the following conditions:\par
\noindent
$(a_1)$\hskip 5pt for each $x\in X$, the function $P(x,\cdot)$ is
quasi-concave and continuous\ ;\par
\noindent
$(a_2)$\hskip 5pt for each $\lambda\in A$, the
function $P(\cdot,\lambda)$ is lower semicontinuous and inf-compact\ ;
\par
\noindent
$(a_3)$ one has
$$\sup_{\lambda\in A}\inf_{x\in X}P(x,\lambda)<
\inf_{x\in X}\sup_{\lambda\in A}P(x,\lambda)\ .$$
\indent
Under such hypotheses, there exists $\lambda^*\in A$ such that the function
$P(\cdot,\lambda^*)$ has at least two global minima.} \par
\medskip
Here is our main result:\par
\medskip
THEOREM 1. - {\it Let $X$ be a reflexive real Banach space;
$I:X\to {\bf R}$ a sequentially weakly lower semicontiunuous and coercive
$C^1$ functional whose derivative admits a continuous inverse on $X^*$;
$J, \Psi, \Phi:X\to {\bf R}$ three $C^1$ functionals with compact derivative
satisfying the following conditions:\par
$$\liminf_{\|x\|\to +\infty}{{J(x)}\over {I(x)}}\geq 0\ ,\hskip 5pt 
\limsup_{\|x\|\to +\infty}{{J(x)}\over {I(x)}}<+\infty\ ,\eqno{(1)}$$
$$\liminf_{\|x\|\to +\infty}{{\Psi(x)}\over {I(x)}}=-\infty\ ,\eqno{(2)}$$
$$\inf_{x\in X}(\Psi(x)+\lambda\Phi(x))>-\infty\eqno{(3)}$$
for all $\lambda>0$. Moreover, assume that there exist a strict local
minimum $x_0$
of $I$, with $I(x_0)=J(x_0)=\Psi(x_0)=\Phi(x_0)=0$, and another point
$x_1\in $X such
that
$$\max\{J(x_1),\Psi(x_1),\Phi(x_1)\}<0\ ,\eqno{(4)}$$
$$\min\left \{ \liminf_{x\to x_0}{{J(x)}\over {I(x)}}, 
 \liminf_{x\to x_0}{{\Phi(x)}\over {I(x)}}\right \} \geq 0\eqno{(5)}$$
and
$$\liminf_{x\to x_0}{{\Psi(x)}\over {I(x)}}>-\infty\ .$$
Under such hypotheses, for each $\nu, \mu$ satisfying
$$\nu>\max\left \{ 0,-{{I(x_1)}\over {J(x_1)}}\right \}\eqno{(6)}$$
and
$$\mu>\max\left \{ 0, -\liminf_{x\to x_0}{{\Psi(x)}\over {I(x)}},
\inf_{r>\sup_{M_{\nu}}\Phi}\mu^*(I+\nu J,\Psi,\Phi,r)\right \}\ ,\eqno{(7)}$$
 where $M_{\nu}$ is the set of all global minima
of $I+\nu J$, there
exists $\lambda^*>0$ such that the functional
$\mu(I+\nu J)+\Psi+\lambda^*\Phi$ has at least four critical points. Precisely,
among them, one is $x_0$ as a strict local, not global minimum and two are global minima.}
\par
\smallskip
PROOF. First of all, observe that,
since $X$ is reflexive, the functionals $J,\Psi,\Phi$ are sequentially weakly
continuous, being with compact derivative ([6], Corollary 41.9).
Fix $\nu$ as in $(6)$. For $x\in X\setminus I^{-1}(0)$,
we have
$$I(x)+\nu J(x)=I(x)\left ( 1+\nu{{J(x)}\over {I(x)}}\right )$$
and so, since $I$ is coercive, in view of $(1)$, it follows that
$$\lim_{\|x\|\to +\infty}(I(x)+\nu J(x))=+\infty\ .\eqno{(8)}$$
By the reflexivity of $X$ again, this implies that the set
$M_{\nu}$ is non-empty and bounded.
 As a consequence, $\Phi$ is bounded
in $M_{\nu}$. Also, by $(2)$ and $(3)$, we have $\sup_X\Phi=+\infty$.
Now, fix $\mu$ as in $(7)$. Let $r>\sup_{M_{\nu}}\Phi$ be such that
$\mu>\mu^*((I+\nu J),\Psi,\Phi,r)$.
Since $\Phi^{-1}(r)$ is non-empty and sequentially weakly
closed, there exists $\bar x\in \Phi^{-1}(r)$ such that
$$I(\bar x)+\nu J(\bar x)=\inf_{x\in \Phi^{-1}(r)}(I(x)+\nu J(x))\ .$$
The choice of $r$ implies that $\bar x\not\in M_{\nu}$. So, we infer that
$$\inf_{x\in \Phi^{-1}(]-\infty,r])}
(I(x)+\nu J(x))<\inf_{x\in \Phi^{-1}(r)}(I(x)+\nu J(x))\ .$$
Moreover,
by $(2)$, there exists a sequence $\{x_n\}$ in $X$ such
that
$$\lim_{n\to \infty}
\|x_n\|=+\infty\ ,\hskip 3pt\lim_{n\to \infty}{{\Psi(x_n)}\over {I(x_n)}}=
-\infty \ .\eqno{(9)}$$
For any $\rho\in {\bf R}$ and for $n$ large enough, we have
$$\rho(I(x_n)+\nu J(x_n))+\Psi(x_n)=
(I(x_n)+\nu J(x_n))\left ( \rho+ {{{{\Psi(x_n)}}\over {I(x_n)}}\over
{{1+\nu {{J(x_n)}\over I(x_n)}}}}\right )\ .\eqno{(10)}$$
Clearly, from $(1)$, $(8)$, $(9)$ and $(10)$, it follows that
$$\lim_{n\to \infty}(\mu(I(x_n)+\nu J(x_n))+\Psi(x_n))=-\infty\ .$$
So, if we consider $X$ endowed with the weak topology, all the assumptions
of Theorem A (with $\Gamma=I+\nu J$) are satisfied, and so we have
$$\alpha(\mu(I+\nu J)+\Psi,\Phi,r)<
\beta(\mu(I+\nu J)+\Psi,\Phi,r)\ .$$
But, by Theorem 1 of [1], this inequality is equivalent to
$$\sup_{\lambda\geq 0}\inf_{x\in X}
((\mu(I(x)+\nu J(x))+\Psi(x)+\lambda(\Phi(x)-r))<
\inf_{x\in X}\sup_{\lambda\geq 0}
((\mu(I(x)+\nu J(x))+\Psi(x)+\lambda(\Phi(x)-r))\ .$$
At this point, after observing that, in view of $(3)$ and $(8)$, one has
$$\lim_{\|x\|\to +\infty}(\mu(I(x)+\nu J(x))+\Psi(x)+\lambda\Phi(x))=+\infty
\eqno{(11)}$$
for all $\lambda>0$, we realize that
 we can apply Theorem B, with
$A=]0,+\infty[$, considering $X$ with the weak topology again and
taking 
$$P(x,\lambda)=\mu(I(x)+\nu J(x))+\Psi(x)+\lambda(\Phi(x)-r))\ .$$
 Therefore, there exists
$\lambda^*>0$ such that the functional $\mu(I+\nu J)+\Psi+\lambda^*\Phi$
 has at
least two gobal minima. 
Now, choose $\epsilon,\sigma>0$                        
so that
$$\liminf_{x\to x_0}{{\Psi(x)}\over {I(x)}}>-\mu+\epsilon$$
and
$$\sigma<{{\epsilon}\over {\mu\nu+\lambda^*}}\ .$$
In view of $(5)$ and recalling that $x_0$ is a strict local minimum of $I$
(with $I(x_0)=0$), we can find a neighbourhood $V$ of $x_0$ such that,
for each $x\in V\setminus \{x_0\}$, one has
$$I(x)>0\ ,$$
$${{\Psi(x)}\over {I(x)}}>-\mu+\epsilon\ ,$$
$${{J(x)}\over {I(x)}}>-\sigma\ ,$$
and
$${{\Phi(x)}\over {I(x)}}>-\sigma\ .$$
Consequently, for each $x\in V\setminus \{x_0\}$, we have
$$\mu(I(x)+\nu J(x))+\Psi(x)+\lambda^*\Phi(x)=
I(x)\left ( \mu+\mu\nu{{J(x)}\over {I(x)}}+{{\Psi(x)}\over {I(x)}}
+\lambda^*{{\Phi(x)}\over {I(x)}}\right ) >$$
$$>I(x)(\epsilon-\sigma(\mu\nu+\lambda^*))>0\ .$$
Hence, since $\mu(I(x_0)+\nu J(x_0))+\Psi(x_0)+\lambda^*\Phi(x_0)=0$, it
follows that $x_0$
is a strict local minimum of the functional
$\mu(I+\nu J)+\Psi+\lambda^*\Phi$. On the other
hand, in view of $(4)$ and $(6)$, we have
$$\mu(I(x_1)+\nu J(x_1))+\Psi(x_1)+\lambda^*\Phi(x_1)<0\ ,$$
and hence $x_0$ is not a global minimum of the functional
$\mu(I+\nu J)+\Psi+\lambda^*\Phi$.
 Now, we remark that
this functional, due to $(11)$ and to our assumptions on $I,J,\Phi, \Psi$
 turns out to satisfy the Palais-Smale condition
([6], Example 38.25). Summarizing: the functional
$\mu(I+\nu J)+\Psi+\lambda^*\Phi$ is $C^1$, satisfies the Palais-Smale
condition, has at least two global minima and admits $x_0$ as
a local, not global minimum. At this point,
we can invoke Theorem (1.ter) of [2]
 to ensure the existence of a fourth critical point
for the same functional, and the proof is complete.
\hfill $\bigtriangleup$\par
\medskip
Now, we are going to present an application of Theorem 1 to
quasilinear elliptic equations.\par
\smallskip
So, let $\Omega\subset {\bf R}^n$ be a bounded domain with smooth
boundary and let $p>1$. On the Sobolev space $W^{1,p}_0(\Omega)$, we
consider the norm
$$\|u\|=\left ( \int_{\Omega}|\nabla u(x)|^p dx\right ) ^{1\over p}\ .$$
If $n\geq p$, we denote by ${\cal A}$ the class of all
Carath\'eodory functions $f:\Omega\times {\bf R}\to {\bf R}$ such that
$$\sup_{(x,\xi)\in \Omega\times {\bf R}}{{|f(x,\xi)|}\over
{1+|\xi|^q}}<+\infty\ ,$$
where  $0<q< {{pn-n+p}\over {n-p}}$ if $p<n$ and $0<q<+\infty$ if
$p=n$. While, when $n<p$, we denote by ${\cal A}$  the class
of all Carath\'eodory functions $f:\Omega\times {\bf R}\to {\bf R}$ such
that, for each $r>0$, the function $x\to \sup_{|\xi|\leq r}|f(x,\xi)|$ belongs
to $L^{1}(\Omega)$.\par
Given $f\in {\cal A}$, consider the following Dirichlet problem
$$\cases {-\hbox {\rm div}(|\nabla u|^{p-2}\nabla u)=
f(x,u)
 & in
$\Omega$\cr & \cr u=0 & on
$\partial \Omega$\ .\cr}\eqno{(P_{f})} $$
 Let us recall
that a weak solution
of $(P_{f})$ is any $u\in W^{1,p}_0(\Omega)$ such that
 $$\int_{\Omega}|\nabla u(x)|^{p-2}\nabla u(x)\nabla v(x)dx
-\int_{\Omega}
f(x,u(x))v(x)dx=0$$
for all $v\in W^{1,p}_0(\Omega)$.\par
\smallskip
The functionals $T, J_{f}:W^{1,p}_0(\Omega)\to {\bf R}$ defined by
$$T(u)={{1}\over {p}}\|u\|^p$$
$$J_{f}(u)=\int_{\Omega}F(x,u(x))dx\ ,$$
where
$$F(x,\xi)=\int_{0}^{\xi}f(x,t)dt\ ,$$
are $C^1$ with derivatives given by
$$T'(u)(v)=\int_{\Omega}|\nabla u(x)|^{p-2}\nabla u(x)\nabla v(x)dx$$
$$J'_{f}(u)(v)=\int_{\Omega}f(x,u(x))v(x)dx$$
for all $u,v\in W^{1,p}_0(\Omega)$. Consequently, the weak solutions 
of problem $(P_{f})$ are exactly the critical points in
$W^{1,p}_0(\Omega)$ of the functional $T-J_{f}$ which is called
the energy functional of problem $(P_f)$. Moreover,
$J'_{f}$ is compact, while $T'$ is a homeomorphism between 
$W^{1,p}_0(\Omega)$ and its dual. \par
\smallskip
The announced application of Theorem 1 is as follows:\par
\medskip
THEOREM 2. - {\it Let $q>p$, with $q<{{pn}\over {n-p}}$ when $n>p$, 
and let
$f, g, h:\Omega\times {\bf R}\to {\bf R}$ be three functions belonging to
${\cal A}$ and satisfying the following conditions:\par
$$\lim_{\xi\to +\infty}{{\inf_{x\in \Omega}F(x,\xi)}\over {\xi^p}}
=+\infty\ ,\hskip 5pt
\limsup_{|\xi|\to +\infty}{{\sup_{x\in \Omega}F(x,\xi)}\over {|\xi|^q}}
<+\infty\ , \eqno{(12)}$$
$$\lim_{|\xi|\to +\infty}{{\inf_{x\in \Omega}G(x,\xi)}\over {|\xi|^q}}
=+\infty\ ,\eqno{(13)}$$
$$\limsup_{|\xi|\to +\infty}{{\sup_{x\in \Omega}H(x,\xi)}\over
{|\xi|^p}}\leq 0\ ,\hskip 5pt
\liminf_{|\xi|\to +\infty}{{\inf_{x\in \Omega}H(x,\xi)}\over {|\xi|^p}}>-\infty
\ ,\eqno{(14)}$$
$$\limsup_{\xi\to 0}{{\sup_{x\in \Omega}F(x,\xi)}\over {|\xi|^p}}<+\infty
\ ,\eqno{(15)}$$
$$\liminf_{\xi\to 0}{{\inf_{x\in \Omega}G(x,\xi)}\over {|\xi|^p}}\geq 0
\ ,\eqno{(16)}$$
$$\limsup_{\xi\to 0}{{\sup_{x\in \Omega}H(x,\xi)}\over {|\xi|^p}}\leq 0\ .
\eqno{(17)}$$
Finally, assume that there exist a measurable set
$B\subset \Omega$, with \hbox {\rm meas}$(B)>0$,
and $\xi_1\in {\bf R}$ such that
$$\max\{-F(x,\xi_1), G(x,\xi_1), -H(x,\xi_1)\}<0$$
for all $x\in B$.\par
Under such hypotheses, for each $\nu>0$ large enough, there exists
$\epsilon_{\nu}>0$ with the following property: for each $\epsilon\in
]0,\epsilon_{\nu}[$ there exists $\lambda^*>0$ such that
the problem
$$\cases {-\hbox {\rm div}(|\nabla u|^{p-2}\nabla u)=\epsilon f(x,u)-
\lambda^* g(x,u)+\nu h(x,u) & in $\Omega$
\cr & \cr
u=0 & on $\partial\Omega$ \cr}$$
has at least three non-zero weak solutions, two of which are global minima
in $W^{1,p}_0(\Omega)$ of the corresponding energy functional.} \par
\smallskip
PROOF. First, observe that from the first assumption in $(12)$ it follows
$$\limsup_{\|u\|\to +\infty}{{J_f(u)}\over {\|u\|^p}}=+\infty\ .$$
This is proved in the proof of Theorem 4 of [5], and so we do not
repeat the argument here. Moreover, from the second assumption in $(12)$
and from $(13)$, it clearly follows that, 
for each $\lambda>0$, the function
$\lambda G-F$ is bounded below in ${\bf R}$ (see [5] again), and so the functional
$\lambda J_g-J_f$ is bounded below in $W^{1,p}_0(\Omega)$. Moreover,
by $(14)$, there is $c>0$ and, for each $\epsilon>0$, another $c_{\epsilon}>0$,
such that
$$-c(|\xi|^p+1)\leq H(x,\xi)\leq\epsilon |\xi|^p+c_{\epsilon}$$
for all $(x,\xi)\in \Omega\times {\bf R}$. This clearly implies that
$$\limsup_{\|u\|\to +\infty}{{J_h(u)}\over {\|u\|^p}}\leq 0$$
and
$$\liminf_{\|u\|\to +\infty}{{J_h(u)}\over {\|u\|^p}}>-\infty\ .$$
Furthermore, by $(15)$ and by the second assumption in $(12)$,
there is a constant $c_1>0$ such that
$$F(x,\xi)\leq c_1(|\xi|^p+|\xi|^q)$$
for all $(x,\xi)\in \Omega\times {\bf R}$. Since $q>p$, this implies that
$$\limsup_{u\to 0}{{J_f(u)}\over {\|u\|^p}}<+\infty\ .$$
Now, suppose $n\geq p$. By $(16)$ and $(17)$, taking into account that
$g, h\in {\cal A}$, for some $r>p$ and
for each $\epsilon>0$ there is $d_{\epsilon}>0$
such that
$$G(x,\xi)\geq -\epsilon |\xi|^p-d_{\epsilon}|\xi|^r$$
and
$$H(x,\xi)\leq \epsilon |\xi|^p+d_{\epsilon}|\xi|^r$$
for all $(x,\xi)\in \Omega\times {\bf R}$.
 From this, we get
$$\liminf_{u\to 0}{{J_g(u)}\over {\|u\|^p}}\geq 0$$
and
$$\limsup_{u\to 0}{{J_h(u)}\over {\|u\|^p}}\leq 0\ .$$
In the case $n<p$, we get again these two inequalities thanks to $(16)$ and
$(17)$ and to the continuous embedding of $W^{1,p}_0(\Omega)$ into
$C^0(\overline \Omega)$. Now, let $\omega\in L^1(\Omega)$ be
such that
$$\max\{|F(x,\xi)|, |G(x,\xi)|, |H(x,\xi)|\}\leq \omega(x)$$
for all $(x,\xi)\in \Omega\times [-|\xi_1|,|\xi_1|]$.
Next, choose
a closed set $C\subset B$, an open set $D\subset \Omega$, with
$C\subset D$, and $\theta\in {\bf R}$ in such a way that
$$\eta:=
\max\left \{ -\int_C F(x,\xi_1)dx, \int_C G(x,\xi_1)dx,
-\int_C H(x,\xi_1)dx\right \}<\theta<0$$
and
$$\int_{D\setminus C}\omega(x)dx<\theta-\eta\ .$$
Finally, let $v_1:\Omega\to [-|\xi_1|,|\xi_1|]$ be a function
belonging to $W^{1,p}_0(\Omega)$ such that $v_1(x)=\xi_1$
for all $x\in C$ and $v_1(x)=0$ for all $x\in \Omega\setminus
D$.
Clearly, we have
$$\max\left \{ -\int_{\Omega} F(x,v_1(x))dx, \int_{\Omega} F(x,v_1(x))dx,
-\int_{\Omega} H(x,v_1(x))dx\right \}<\eta+(\theta-\eta)<0\ .$$
At this point,
the conclusion comes directly from that of Theorem 1 applied taking
$I(u)={{1}\over {p}}\|u\|^p$, $J(u)=-J_h(u)$,
$\Psi(u)=-J_f(u)$,
$\Phi(u)=J_g(u)$ for all $u\in W^{1,p}_0(\Omega)$.\hfill $\bigtriangleup$
\bigskip
\bigskip
\centerline {\bf References}\par
\bigskip
\bigskip
\noindent
[1]\hskip 5pt G. CORDARO, {\it On a minimax problem of Ricceri},
J. Inequal. Appl., {\bf 6} (2001), 261-285.\par
\smallskip
\noindent
[2]\hskip 5pt N. GHOUSSOUB and D. PREISS, {\it A general mountain pass principle
for locating and classifying critical points}, Ann. Inst. H. Poincar\'e
Anal. Non Lin\'eaire, {\bf 6} (1989), 321-330.\par
\smallskip
\noindent
[3]\hskip 5pt B. RICCERI, {\it Multiplicity of global minima for
parametrized functions}, Rend. Lincei Mat. Appl., {\bf 21} (2010),
47-57.\par
\smallskip
\noindent
[4]\hskip 5pt  B. RICCERI, {\it A class of nonlinear eigenvalue problems with
four solutions}, J. Nonlinear Convex Anal., {\bf 11} (2010), 503-511.\par
\smallskip
\noindent
[5]\hskip 5pt B. RICCERI, {\it A further refinement of a three critical points
theorem}, Nonlinear Anal., {\bf 74} (2011), 7446-7454.\par
\smallskip
\noindent
[6]\hskip 5pt E. ZEIDLER, {\it Nonlinear functional analysis and its
applications}, vol. III, Springer-Verlag, 1985.\par
\bigskip
\bigskip
\bigskip
\bigskip
Department of Mathematics\par
University of Catania\par
Viale A. Doria 6\par
95125 Catania\par
Italy\par
{\it e-mail address}: ricceri@dmi.unict.it

\bye